\documentclass[10pt]{amsart}
\usepackage{amssymb}
\usepackage{MnSymbol}

\usepackage{amscd}
\usepackage{amsthm}
\usepackage[dvips]{graphicx}
\usepackage{color}
\usepackage[all]{xy}
\usepackage{verbatim}

\newtheorem{Lemma1}{{Lemma}}[section]
\newtheorem{Theo1}[Lemma1]{{Theorem}}
\newtheorem*{Theo2}{{Theorem}}
\newtheorem{Def1}[Lemma1]{{Definition}}
\newtheorem{Prop1}[Lemma1]{{Proposition}}
\newtheorem{Claim1}[Lemma1]{{Claim}}
\newtheorem{Rem1}[Lemma1]{{Remark}}
\newtheorem{Cor1}[Lemma1]{{Corollary}}
\newtheorem{Ex1}[Lemma1]{{Example}}
\newtheorem{Not1}[Lemma1]{{Notation}}

\newenvironment{Lemma}{\begin{Lemma1}}{\end{Lemma1}}
\newenvironment{Def}{\begin{Def1}\rm}{\end{Def1}}
\newenvironment{Prop}{\begin{Prop1}}{\end{Prop1}}
\newenvironment{Rem}{\begin{Rem1}\rm}{\end{Rem1}}
\newenvironment{Theorem}{\begin{Theo1}}{\end{Theo1}}

\newenvironment{Cor}{\begin{Cor1}}{\end{Cor1}}

\newenvironment{Example}{\begin{Ex1}\rm}{\end{Ex1}}

\title{Differential graded division algebras, their modules, and dg-simple algebras}

\author{Alexander Zimmermann}
\address{\newline
Universit\'e de Picardie,
\newline D\'epartement de Math\'ematiques et LAMFA (UMR 7352 du CNRS),
\newline 33 rue St Leu,
\newline F-80039 Amiens Cedex 1,
\newline France}
\email{alexander.zimmermann@u-picardie.fr}

\date{August 8, 2024}

\newcommand{\lra}{\longrightarrow}

\newcommand{\ra}{\rightarrow}
\newcommand{\sdp}{\times\kern-.2em\vrule height1.1ex depth-.05ex}
\newcommand{\epi}{\lra \kern-.8em\ra}

\newcommand{\N}{{\mathbb N}}

\newcommand{\Z}{{\mathbb Z}}

\newcommand{\im}{\textup{im}}
\newcommand{\dickebox}{{\vrule height5pt width5pt depth0pt}}

\newcommand{\id}{\textup{id}}

\newcommand{\Hom}{\textup{Hom}}

\newcommand{\End}{\textup{End}}

 \normalbaselines

\begin{document}

\begin{abstract}
We give the definition of a dg-division algebra, that is a concept of a differential graded algebra which
may serve as an analogue of a division algebra. We classify them completely, and 
show that they are either acyclic or have differential $0$.
Further, we prove that the graded centre of dg-simple dg-algebras is a dg-division algebra, and also
the dg-endomorphism ring of a dg-simple module is a dg-division algebra.
We also shall give a Jacobson-Chevalley density theorem for acyclic dg-algebras.
\end{abstract}

\maketitle

\section*{Introduction}

Differential graded algebras (dg-algebras for short) and their differential graded 
modules appear in various places, 
mainly of geometric and topological nature. More precisely, let $K$ be a commutative ring. 
A dg-algebra over $K$ is a $\Z$-graded algebra with a graded $K$-linear 
endomorphism $d$ of degree $1$ with $d^2=0$ satisfying the Leibniz equation 
$d(ab)=d(a)\cdot b+(-1)^{|a|}a\cdot d(b)$ for all homogeneous $a,b\in A$. 
Dg-modules are defined similarly. For more precise definitions we refer 
to Section~\ref{Definitionsection}. 
Though defined by Cartan \cite{Cartandg} in 
1954 already, the ring theory of dg-algebras remained largely unexplored until 
quite recently. In 2002 Aldrich and Garcia Rozas characterized the dg-algebras 
whose dg-module category is semisimple. They obtained that this is the case 
if and only if the algebra is acyclic and the algebra of cycles $\ker(d)$ is 
semisimple as graded modules. Orlov studied in \cite{Orlov1, Orlov2} finite dimensional 
dg-algebras over a field $K$ with a geometric motivation. Using this 
approach Goodbody \cite{Goodbody} studied a version of a dg-Jacobson radical.
In a sequel of papers we studied more systematically the ring theory 
of dg-algebras. In \cite{dgorders} we defined dg-Jacobson radicals in a 
much more general and natural setting,  study dg-simple dg-algebras, dg-simple dg-modules,
and the relation with dg-simple dg-algebras. We study dg-orders and define
locally free dg-class groups. In \cite{dgGoldie} we define Ore localisation
of dg-algebras, proving that under some hypotheses, the localised ring is again dg.
Using this, we study and define dg-uniform dimension, and a dg-Goldie theorem 
for dg-Goldie rings in \cite{dgGoldie}. 

In the present paper we first define a dg-division algebra as a dg-ring without non trivial left 
or right dg-ideals. The question of an appropriate concept for a dg-division algebra was posed 
recently by Violeta Borges Marques and Julie Symons.
We show in Theorem~\ref{dgfieldscharacteristation} 
that under a certain  mild technical condition dg-division algebras are precisely those 
for which the ring of cycles is a $\Z$-graded-division algebras
(cf \cite[page 38]{gradedrings} for the graded concept). For this result we need a technical 
assumption, namely that the set of left regular elements of the ring of cycles coincides with 
the set of right regular elements. This holds true, by a result of Goodearl 
and Stafford \cite{GoodearlStafford} for graded-Noetherian  graded-prime rings of cycles,
using a graded version of Goldie's result.  As a second  main result we 
then classify completely dg-division rings 
under these conditions in Theorem~\ref{acyclicordequalzero}. Namely, they are
either acyclic or else they have differential $0$. The structure is then 
completely given by a result from 
Aldrich and Garcia Rozas \cite{Tempest-Garcia-Rochas} and Nastacescu and van 
Oystaen \cite[A.1.4.3 Corollary]{gradedrings}.

Further, we prove that endomorphism complexes of dg-simple dg-modules over dg-algebras 
are dg-division algebras. We provide further examples. 
We show that  the graded centre of a dg-simple algebra
is a dg-division algebra.

In a final section we show a dg-version of the 
Jordan-Chevalley density theorem in case of acyclic dg-algebras. 

The paper is organized as follows. In Section~\ref{Elementarydefsect} we recall the definition
of dg-algebras and their modules,  the constructions we use in the paper, and we define dg-division algebra, 
the main 
structure studied in our paper.  In Section~\ref{dgfieldsandgrfields} we study the main 
properties of dg-division algebras with respect to the cycles, and some occurrences of dg-division algebras. 
In Section~\ref{acyclicsection} we give additional properties of acyclic dg-algebras, which is then used in 
Section~\ref{dgJacobsonsect} where w prove a version of the Jacobson-Chevalley 
density theorem for acyclic dg-algebras.

\subsubsection*{Acknowledgement:}
I wish to thank Violeta Borges Marques and Julie Symons for having asked me during  
the Oberwolfach workshop 'Hochschild (co-)homology and applications' in April 2024 for
a possible concept  of a dg-field.  

\section{Elementary Definitions}

\label{Elementarydefsect}

\subsection{Differential graded algebra and their modules; definitions and notations}

\label{Definitionsection}

In this subsection we recall notations and basic conventions, as well as the 
construction of a tensor product of two dg-algebras over a common graded commutative 
subalgebra of their graded centres. This should be well-known, but we could not find a reference,
and in any case the reader may appreciate an explicit verification right in the paper.
 
Let $K$ be a commutative ring. A differential graded $K$-algebra (dg-$K$-algebra for short) 
$(A,d)$ is a $\Z$-graded 
algebra $A$ with a graded $K$-linear endomorphism of degree $1$ with $d^2=0$
and satisfying the Leibniz rule $d(a\cdot b)=d(a)\cdot b+(-1)^{|a|}a\cdot d(b)$. 
By definition a dg-ring is a dg-$\Z$-algebra. 

If $(A,d)$ is a dg-$K$-algebra, then define $A^{op}$ to be the same graded 
additive group as $A$, but 
multiplication defined as $a\cdot_{op}b:=(-1)^{|a|\cdot|b|}b\cdot a$. 
Then $(A^{op},d)$ is again a dg-$K$-algebra. 

A differential graded left module $(M,\delta)$ over $(A,d)$ is a $\Z$-graded 
$A$-module $M$ with an endomorphism $\delta$ of degree $1$ 
satisfying $\delta(a\cdot m)=d(a)\cdot m+(-1)^{|a|}a\cdot \delta(m)$
for all homogeneous $a\in A$ and $m\in M$. Occasionally we 
denote a differential graded module by dg-module for short. 

A differential graded right module over $(A,d)$ is a differential graded left module over $(A^{op},d)$. 

For two differential graded left modules $(M,\delta_M)$ and $(N,\delta_N)$ over $(A,d)$ 
we put 
$$\Hom_A^\bullet((M,\delta_M),(N,\delta_N)):=\{f\in \Hom_{graded}(M,N)\;|\;f(am)=(-1)^{|a|\cdot|f|}af(m)\}$$
and $$d_\Hom(f):=\delta_N\circ f-(-1)^{|f|}f\circ\delta_M.$$
Further, $\End_A^\bullet(M,\delta_M):=\Hom_A^\bullet((M,\delta_M),(M,\delta_M)).$
Then $(\End_A^\bullet(M,\delta_M),d_\Hom)$ is a dg-$K$-algebra, and
 $(\Hom_A^\bullet((M,\delta_M),(N,\delta_N)),d_\Hom)$
is a left dg-module over $\End_A^\bullet(N,\delta_N)$ and a right dg-module over 
$\End_A^\bullet(M,\delta_M)$.

\medskip

We further mention that for a graded algebra $A$ by some group $G$
we denote be a gr-simple $A$-module an $A$-module which 
does not allow a $G$-graded submodule other than $0$ or itself. Similarly, we denote 
by  a gr-simple algebra an algebra which is gr-simple as $A\otimes A^{op}$-module, 
and a gr-division algebra an algebra where all homogeneous non zero elements are invertible. 
Occasionally we use the notion graded-simple for gr-simple, etc.

\subsection{Elements of the ring theory of differential graded algebra}

\begin{Def}
A dg-$K$-algebra $(A,d)$ is {\em dg-simple} if $(A,d)$ does not contain any twosided dg-ideal 
other than $0$ and $A$.
\end{Def}

\begin{Def}
A dg-$K$-algebra $(A,d)$ is a {\em dg-division algebra} if $(A,d)$ does not contain any dg-left 
ideal nor a dg-right ideal other than $0$ and $A$. 
\end{Def}

\begin{Rem}
In case we want to stress that a dg-division algebra is
commutative (or graded commutative) we shall call it a  dg-field (or 
graded commutative dg-field). 
\end{Rem}


\begin{Example}
Let $K$ ba a field. Then by \cite{dgGoldie} 
the dg-ring $K[X,X^{-1}]$ with differential $d(X^{2n+1})=X^{2n}$ and $d(X^{2n})=0$ for all $n\in\Z$ 
is dg-simple. Since it is commutative, it is a dg-field.
Note however that the ungraded ring $K[X,X^{-1}]$ is not Artinian (though Noetherian). 
\end{Example}

Acyclic dg-algebras will play a crucial r\^ole. They are classified in \cite{Tempest-Garcia-Rochas}.

\begin{Theorem} \cite[Theorem 4.7; Proposition 4.3]{Tempest-Garcia-Rochas}
\label{acyclicalgebras}
Let $(A,d)$ be a dg-algebra over a commutative ring $R$. Then the following statements are equivalent:
\begin{enumerate}
\item $(A,d)$ is acyclic
\item The left regular module is a projective object in the category of left dg-modules
\item $1\in\im(d)$
\item Taking cycles is a right exact functor $Z(-)$ from the category of left $(A,d)$-dg-modules to the category of graded modules over $\ker(d)$.
\item \label{acyclicalgebras5}$(A\otimes_{\ker(d)}-Z(-))$ are mutually inverse equivalences between 
the category of left $(A,d)$-dg-modules to the category of graded modules over $\ker(d)$.
\end{enumerate}
Further, if $(A,d)$ is acyclic, and $d(y)=1$ for some homogeneous element $y\in A$, then 
we have an isomorphism of dg-algebras between $A$ and the following quotient of the twisted 
polynomial algebra
$$A\simeq\ker(d)[X;D]/(X^2-y^2)$$
with $D(a):=(-1)^{|a|}d(yay)$, for all homogeneous $a\in\ker(d)$. 
\end{Theorem}

\section{Relating dg-division algebras with gr-division algebra}

\label{dgfieldsandgrfields}

\begin{Lemma}\label{cyclesarefieldimpliesdgfield}
Let $(A,d)$ be a dg-ring  and suppose that $\ker(d)$ is a $\Z$-gr division algebra. Then $(A,d)$ is a 
dg-division algebra. 
\end{Lemma}

Proof. Let $(I,d)$ be a non trivial  dg-left ideal of $(A,d)$. 
Then $I\cap\ker(d)$ is a $\Z$-graded ideal of $\ker(d)$.

Since $I\neq 0$, there is $0\neq x\in I$. If $x\in\ker(d)$, we have shown
that $I\cap\ker(d)\neq 0$. If $d(x)\neq 0$, then $d(x)\in I$ since
$I$ is a dg-ideal, and $d(x)\in\ker(d)$ since $d^2=0$. Hence $I\cap\ker(d)\neq 0$. 

Suppose that  $I\cap\ker(d)=\ker(d)$. Since $d(1)=0$, we get $1\in I\cap\ker(d)$. Hence $1\in I$
and therefore $A=I$. This is a contradiction to the hypothesis that $I$ is a non trivial dg-ideal. 
The same holds for dg-right ideals.
\dickebox

\begin{Lemma}\label{Adgfieldthencyclesgrfield}
Let $(A,d)$ be a dg-ring and suppose that $(A,d)$ does not contain any non trivial dg-left ideal. 
If the set of homogeneous right regular  element of $\ker(d)$ coincides 
with the set of homogeneous left regular elements of $\ker(d)$, then 
$\ker(d)$ is a $\Z$-gr-division algebra.  
\end{Lemma}

Proof. Let $u\neq 0$ be a homogenous element of $\ker(d)$. Then $Au$ is a dg-left ideal of $A$. 
Since $u\neq 0$, we get $Au\neq 0$.
Suppose that $Au=A$. Then there is $b\in A$ with $1=bu.$
Then $$0=d(1)=d(b)\cdot u+(-1)^{|b|}b\cdot d(u)=d(b)\cdot u.$$
Since $u$ has a left inverse in $A$, it is right regular. Indeed, if $ux=0$ in $A$, then $x=bux=0$.
Hence, by hypothesis, $u$ is left regular as well. Therefore, $d(b)=0$ and hence $b\in\ker(d)$.
This shows that the equation $bu=1$ already holds in $\ker(d)$, and hence $u$ is left invertible
in $\ker(d)$. Therefore, $\ker(d)$ does not contain any non trivial graded left ideal, and by
\cite[page 38; Lemma 1.4.1]{gradedrings}  
$u$ is invertible in $\ker(d)$. \dickebox

\begin{Rem} \label{leftregularequalsrightregular}
 \label{whenthestrangeconditionholds}
We may consider situations where the hypothesis that
the set of homogeneous right regular  element of $\ker(d)$ coincides 
with the set of homogeneous left regular elements of $\ker(d)$ holds. 
\begin{itemize}
\item 
If $\ker(d)$ is finite dimensional over some field, this is true, since multiplication by a
left regular element $u$ is realised by a matrix with non zero determinant.  
\item
If $\ker(d)$ is graded commutative, or commutative, then trivially this property holds. 
\item
It is known that if $\ker(d)$ is left gr-Noetherian and gr-prime, this is true. Indeed, by 
\cite{GoodearlStafford} localising at the homogeneous regular elements of $\ker(d)$ 
yields a gr-simple algebra with all homogeneous regular elements being invertible. 
Now, the proof of \cite[(5.8) Proposition, (5.9) Proposition]{CR1} applies verbatim to graded 
rings and homogeneous left invertible elements. 
\end{itemize}
\end{Rem}

\begin{Theorem} \label{dgfieldscharacteristation}
Let $(A,d)$ be a differential graded ring. Suppose that in $\ker(d)$ the set of the left 
regular homogeneous elements and the set of the right regular homogeneous elements coincide.
Then $(A,d)$ is a dg-division algebra if and only if $\ker(d)$ is a gr-division algebra. 
\end{Theorem}

Proof. This is a direct consequence of Lemma~\ref{cyclesarefieldimpliesdgfield} and 
Lemma~\ref{Adgfieldthencyclesgrfield}. Note that the lemmas give a more general statement for 
the only if direction. \dickebox

\begin{Cor} \label{tripletsequiv}
Let $(A,d)$ be a differential graded ring. Suppose that in $\ker(d)$ the set of the left 
regular homogeneous elements and the set of the right regular homogeneous elements coincide.
Then the following statements are equivalent. 
\begin{itemize}
\item  $(A,d)$ admits only trivial dg-left ideals.
\item  $(A,d)$ admits only trivial dg-right ideals.
\item each non zero homogeneous element of $\ker(d)$ is invertible (i.e. $\ker(d)$ is a gr-division algebra).
\end{itemize} 
\end{Cor}

Proof. The first statement implies the third statement by Lemma~\ref{Adgfieldthencyclesgrfield}. 
Similarly, the second statement implies the third statement by Lemma~\ref{Adgfieldthencyclesgrfield}. 
The third statement implies each of the first and the second 
statement by Lemma~\ref{cyclesarefieldimpliesdgfield}.
\dickebox

\medskip

We shall need to study the graded centre of a dg-algebra. Recall the definition. 
Let $A$ be a $\Z$-graded algebra. Then the graded centre $Z_{gr}(A)$ is defined as
$$Z_{gr}(A):=\{b\in A\;|\;b\textup{ homogeneous and }ba=(-1)^{|a||b|}ab\textup{ for all homogeneous }b\in A\}$$

\begin{Lemma}\label{Zgrisadgalg}
Let $(A,d)$ be a differential graded ring. Then $(Z_{gr}(A),d)$ is a differential graded subalgebra.
\end{Lemma}

Proof. By construction, $Z_{gr}(A)$ is a subalgebra.
Further, for any $a\in Z_{gr}(A)$  we get
\begin{eqnarray*}
d(a)b&=&d(ab)-(-1)^{|a|}a\cdot d(b)\\
&=&(-1)^{|a||b|}d(ba)-(-1)^{|a|+(|a|(|b|+1))}d(b) a\\
&=&(-1)^{|a||b|}\left(d(ba)-d(b) a\right)\\
&=&(-1)^{|a||b|+|b|}b\cdot d(a)\\
&=&(-1)^{|d(a)||b|}b\cdot d(a)
\end{eqnarray*}
Hence $(Z_{gr}(A),d|_{Z_{gr}(A)})$ is a dg-algebra again.
\dickebox

\begin{Rem} \label{gradedcommutativeinoddchar}
Recall that in a graded commutative ring, the square of homogeneous elements of odd 
degree is $0$, unless the characteristic of the base ring is $2$. Indeed, let $x$ be 
a homogeneous element of odd degree. Then,
swapping the two factors yields 
$$x\cdot x=(-1)^{|x|\cdot|x|}x\cdot x=-x\cdot x.$$
Hence $2\cdot x^2=0$. 
\end{Rem}

\begin{Cor} \label{gradedcenterofdgfieldisdgfield}
Let $(D,\partial)$ be a dg-division algebra. Suppose that in $\ker(\partial)$ the set of the left 
regular homogeneous elements and the set of the right regular homogeneous elements coincide. 
Then $(Z_{gr}(D),\partial)$ is a dg-division algebra as well. 

If $(D,\partial)$ is a dg-algebra, and if a homogeneous 
$z\in Z_{gr}(D)$ has a homogeneous inverse in $\ker(\partial)\subseteq A$, then $a\in Z_{gr}(D)$.

In particular, if the characteristic of $D$ is 
different from $2$, then $Z_{gr}(D)\cap\ker(\partial)$ is concentrated in even degrees, and in any 
case $Z_{gr}(D)\cap\ker(\partial)$ is commutative and either isomorphic to a Laurent polynomial ring
or concentrated in degree $0$. 
\end{Cor}

Proof. 
By Lemma~\ref{Zgrisadgalg} we see that $Z_{gr}(D)$ is a dg-subalgebra of $(D,\partial)$. 
Further, all elements in $Z_{gr}(D)\cap\ker(\partial)$ are invertible in $\ker(d)$. 
However, if $z\in Z_{gr}(D)$ and $a\in \ker(\partial)$ are homogeneous elements 
with $az=za=1$ (using that in a group the left and the right inverse of a fixed element coincide), 
then for all homogeneous $b\in D$
$$z(ab-(-1)^{|a||b|}ba)=zab-(-1)^{|a||b|}zba=b-(-1)^{|a||b|+|z||b|}bza=b-b=0$$
using that $|a|=-|z|$. Hence $a\in Z_{gr}(D)$ as well. 
Using Corollary~\ref{tripletsequiv} and Remark~\ref{gradedcommutativeinoddchar} 
we proved the statement. The classification of commutative $\Z$-graded gr-division algebras 
is given in \cite[Section A.1.4; Section B.1.1]{gradedrings}, and it is shown that 
 $Z_{gr}(D)\cap\ker(\partial)$ is isomorphic to a Laurent polynomial ring or concentrated in 
 degree $0$.
\dickebox

\begin{Theorem}\label{acyclicordequalzero}
Let $(A,d)$ be a dg-division algebra and suppose that the set of homogeneous left regular elements 
of $\ker(d)$ coincides with the set of homogeneous right regular elements of $\ker(d)$. Then either 
$(A,d)$ is acyclic or $d=0$. Further, $\ker(d)=:R$ is either concentrated in degree $0$, or else
$\ker(d)$ is a twisted Laurent polynomial ring 
$\ker(d)=R_0[T,T^{-1};\sigma]$ for 
some automorphism $\sigma$ of $R_0$ 
with $Xr=\sigma(r)X$ for all $r\in R_0$.
\end{Theorem}

Proof.
Indeed, since $H(A,d)$ is a dg-algebra and since 
$\ker(d)/\im(d)=H(A,d)$ as an algebra, $\im(d)$ is a twosided ideal of $\ker(d)$. 
However, as $\ker(d)$ is a gr-division algebra, every homogeneous non zero element of $\ker(d)$ is invertible in $\ker(d)$. Therefore, either $\ker(d)=\im(d)$, or
$\im(d)$ does not contain any non zero homogeneous element of $\ker(d)$. Hence, either $H(A,d)=0$ in 
the first case, or in the second case $d=0$. The structure of $\ker(d)$ follows from Nastacescu-van Oystaen~\cite[A.1.4.3 Corollary]{gradedrings}.
\dickebox

\begin{Example}
Both cases mentioned in Theorem \ref{acyclicordequalzero} do occur. 
Let $A=K[X,X^{-1}]$ the Laurent polynomial ring over a field $K$. 
Then the differential $0$ yields a gr-field, whence also a dg-field.
However, also $d(X^{2n+1}):=X^{2n}$ and $d(X^{2n})=0$ for all $n\in\N$
yields a dg-algebra, which is acyclic, and a dg-field, since it is already a gr-field. But also 
the cycles is a Laurent polynomial algebra over $K$ in all even degrees, which is again a gr-field.  
\end{Example}

\begin{Cor} \label{homologyofdgfields}
Let $(A,d)$ be a dg-division algebra, and suppose that the set of homogeneous left regular elements 
of $\ker(d)$ coincides with the set of homogeneous right regular elements of $\ker(d)$. 
Then, either $(A,d)$ is acyclic, or $\ker(d)=H(A,d)$ is a gr-division algebra. 
\end{Cor}

Proof. By Theorem~\ref{dgfieldscharacteristation} 
we see that $\ker(d)$ is a gr-division algebra, whence any homogeneous non zero element is invertible. 
Therefore, also in $H(A,d)=\ker(d)/\im(d)$ every homogeneous non zero element is invertible.
This shows the statement. \dickebox

\medskip

Further, we recall

\begin{Lemma} \cite[Lemma 4.1, Lemma 4.2]{Tempest-Garcia-Rochas} \label{characterizeacylic}
Let $(A,d)$ be a dg-algebra. Then $(A,d)$ is acyclic if and only if $1\in\im(d)$. 

Let $(A,d)$ be an acyclic dg-algebra over a commutative ring $R$. Then, as an $R$-module
$A=\ker(d)\oplus\ker(d)\cdot y=\ker(d)\oplus y\cdot\ker(d)$ for any homogeneous 
element $y\in A$ with $d(y)=1$. 
\end{Lemma}

\begin{Rem}
Note that using Lemma~\ref{characterizeacylic}, if $(A,d)$ is acyclic, then 
 \cite[Proposition 4.3]{Tempest-Garcia-Rochas} gives more explicitly
$A\simeq\ker(d)[X;D]/(X^2-y^2)$, with $d(y)=1$ and $D(a)=-(-1)^{|a|}d(yay)$ 
for any homogeneous $a\in\ker(d)$. Furthermore, 
$A=\ker(d)\oplus\ker(d)y$, and the isomorphism 
is $\Phi:\ker(d)[X;D]\lra A$ by $\Phi(X)=y$ and $\Phi(a)=a$ for any $a\in\ker(d)$.
For any homogeneous $a,b\in\ker(d)$ we get $d(b+ay)=a$. Moreover,
$D(a)=Xa-(-1)^{|a|}aX$ for any homogeneous $a\in\ker(d)$. 
Recall that a dg-division algebra has either differential $0$ or $(A,d)$ is acyclic.
The cycles are still given by $\ker(d)$, which need to be a gr-division ring. 
\end{Rem}

\begin{Prop} \label{endosofsimpleisdgfield}
Let $K$ be a field, 
let $(A,d)$ be a differential graded algebra and let $(S,\delta)$ be a 
dg-simple left dg-module over $(A,d)$. 
Then  $(End_A^\bullet(S,\delta),d_\textup{Hom})$ is a dg-division algebra.
Moreover, the set of homogeneous left regular elements of $\ker(d_\Hom)$ 
coincides with the set of homogeneous right regular elements of $\ker(d_\Hom)$. 
\end{Prop}

Proof.
We shall use Corollary~\ref{tripletsequiv} and we shall see that the hypotheses are satisfied.. 

Now, let $f\neq 0$ be homogeneous with 
$f\in\ker(d_\Hom)$. Then this is equivalent with 
$$0=d_\Hom(f)=\delta\circ f-(-1)^{|f|}f\circ\delta$$
and hence $$f\circ\delta=(-1)^{|f|}\delta\circ f.$$
We claim that $\ker(f)$ is a dg-submodule of $(S,\delta)$.
Let $x\in\ker(f)$ be homogeneous.  
$$
0=\delta(f(x))=(-1)^{|f|} f(\delta(x))
$$
and hence $\delta(x)\in\ker(f)$. 
Moreover, for any homogeneous $a\in A$
 $$f(ax)=(-1)^{|a|\cdot|f|}a\cdot f(x)=0$$
and hence $ax\in\ker(f)$ again. 

Since $(S,\delta)$ is dg-simple, either $f=0$ or $\ker(f)=0$.

Also $\im(f)$ is a dg-submodule of $(S,\delta)$. If $x=f(y)$, then 
$$\delta(x)=\delta(f(y))=(-1)^{|f|}f(\delta(y))\in\im(f)$$
Further, 
$$a\cdot f(x)=(-1)^{|a||f|}f(a\cdot x)\in\im(f)$$
again. Since $(S,\delta)$ is dg-simple, and since $\im(f)$ is a 
dg-submodule, either $\im(f)=0$ (which is equivalent with $f=0$ 
and this was excluded) or $\im(f)=S$. Therefore $f$ is surjective 
as well. 

Hence $f$ is an isomorphism between $(S,\delta)$ and a shifted copy, whence invertible,  
and this shows that $\ker(d_\Hom)$ is a gr-division algebra. 

We want to apply Corollary~\ref{tripletsequiv}. If $f$ is a left invertible 
non invertible dg-endomorphism of $(S,\delta)$, then $(S,\delta)$ has a 
non trivial direct factor, which is contradictory to $(S,\delta)$ being simple. 
Likewise a right invertible non invertible dg-endomorphism leads to a contradiction. 

Hence $(End_A^\bullet(S,\delta),d_\textup{Hom})$ is a dg-division algebra
by Corollary~\ref{tripletsequiv}.
\dickebox

\medskip

We summarize the main result of this section.

\begin{Theorem} \label{dgdivisionalgebraclassification}
Let $(A,d)$ be a dg-algebra. Suppose that the set of left regular homogeneous 
elements of $\ker(d)$ coincides with the set of right regular homogeneous elements of $\ker(d)$.

Then 
\begin{itemize}
\item
$(A,d)$ is a dg-division algebra if and only if $\ker(d)$ is a $\Z$-gr-division algebra (cf \cite{gradedrings}).
\item 
If $(A,d)$ is a dg-division algebra,
\begin{itemize} 
\item then,
\begin{itemize}
\item  either $d=0$ and $\ker(d)$ is  a skew field 
concentrated in degree $0$, 
\item or $H(A,d)=0$ and there is a skew field $R_0$
such that $\ker(d)\simeq R_0[X,X^{-1};\phi]$ for an automorphism $\phi$ of $R_0$ 
and $Xr=\phi(r)X$ for any $r\in R_0$. 
\end{itemize}
\item If $H(A,d)=0$, then also $A\simeq\ker(d)[T;D]/(T^2-y^2)$, with $d(y)=1$ and $D(a)=-(-1)^{|a|}d(yay)$ 
for any homogeneous $a\in\ker(d)$. Furthermore, 
$A=\ker(d)\oplus\ker(d)y$, and the isomorphism 
is $\Phi:\ker(d)[T;D]\lra A$ by $\Phi(T)=y$ and $\Phi(a)=a$ for any $a\in\ker(d)$.
For any homogeneous $a,b\in\ker(d)$ we get $d(b+ay)=a$. Moreover,
$D(a)=Ta-(-1)^{|a|}aT$ for any homogeneous $a\in\ker(d)$.
\end{itemize}
\end{itemize}
\end{Theorem}

\begin{Cor}
Let $(A,d)$ be a dg-algebra. Suppose that the set of left regular homogeneous 
elements of $\ker(d)$ coincides with the set of right regular homogeneous elements of $\ker(d)$.
Suppose that $(A,d)$ is a graded commutative dg-division algebra. 
Then either $\ker(d)$ is a field concentrated in degree $0$, or else 
$\ker(d)=K[X,X^{-1}]$ for some field $K$
and, in case $K$ is of characteristic different from $2$, then $X$ is in even degree.
\end{Cor}

Indeed, since $(A,d)$ is a dg-division algebra, we get that 
either $\ker(d)$ is concentrated in degree $0$, or there is a skew field $K$
such that $\ker(d)\simeq K[X,X^{-1};\phi]$ for an automorphism $\phi$ of $K$ 
and $Xr=\phi(r)X$ for any $r\in K$. 
Now, since $r\in K$ is of degree $0$, and since $A$ is graded commutative,
$$\phi(r)X=Xr=(-1)^{|r|\cdot |X|}rX=rX$$
and hence $\Phi(r)=r$ for all $r\in K$.
Since $A$ is graded commutative, $$X\cdot X=(-1)^{|X|\cdot|X|}X\cdot X$$ 
interchanging the two factors, and hence $|X|$ has to be even. 
\dickebox



\begin{Lemma}\label{gradedcenterisdgfield}
Let $(A,d)$ be a dg-simple dg-algebra. Then $Z_{gr}(A,d)$ is a dg-division algebra. 
\end{Lemma}

Proof.
Let $u\in\ker(d)\cap Z_{gr}(A)$ be a homogeneous non zero element. Then $Au$ is a 
twosided dg-ideal of $A$, using that $u$ is in the graded centre. 
Since $(A,d)$ is dg-simple, there is a homogeneous $a\in A$ with $au=1$. Then 
$$0=d(1)=d(au)=d(a)u+(-1)^{|a|}a d(u)=d(a)u.$$ 
Since $au=1$, and since $u$ is in the graded centre, also $ua=\pm 1$, and 
hence multiplying by $a$ from the right, we get $a\in\ker(d)$. 
Further, $a\in Z_{gr}(A)$, using  Corollary~\ref{gradedcenterofdgfieldisdgfield}. 
Hence, by Lemma~\ref{cyclesarefieldimpliesdgfield} we get that $Z_{gr}(A)$ is 
a dg-division algebra. \dickebox

\medskip

Let $(A,d)$ be a dg-$K$-algebra. Then the graded centre
$Z_{gr}(A,d)$ is a dg-subalgebra 
of $(A,d)$ (cf Lemma~\ref{Zgrisadgalg}). 

If $(A,d_A)$ and $(B,d_B)$ are dg-$K$-algebras, then $((A\otimes_KB),d_{A\otimes_KB})$
is a dg-$K$-algebra again with 
$$d_{A\otimes B}=d_A\otimes \textup{id}_B+\id_A\otimes d_B$$
and 
$$(a_1\otimes b_1)\cdot (a_2\otimes b_2)=(-1)^{|b_1||a_2|}(a_1a_2\otimes b_1b_2).$$

\begin{Lemma}
Let $(A,d_A)$ and $(B,d_B)$ be dg-$K$-algebras, and let $Z=(Z,d)$ be a common graded commutative 
dg-subalgebra of the graded centre of $A$ and of $B$. 
Then there is a differential $d_{A\otimes_ZB}$ on $A\otimes_ZB$ such that $(A\otimes_ZB,d_{A\otimes_ZB})$ 
is a dg-algebra again. 
\end{Lemma}

Proof. First, we consider the ordinary tensor product $A\otimes_ZB$. 
This is graded by posing the degree $n$ component as 
$\bigoplus_{k\in\Z}A_k\otimes_ZB_{n-k}$. This is well-defined since $Z$ is graded. 
Indeed, if $a\in A_k$ is in degree $k-\ell$, $b\in B_{n-k}$ is in degree $n-k$, 
and $z\in Z_\ell$ is in degree $\ell$, then $az\otimes b=a\otimes zb$ and the 
term on both sides is in degree $n$.
Further, put again 
$$(a_1\otimes b_1)\cdot (a_2\otimes b_2)=(-1)^{|b_1||a_2|}(a_1a_2\otimes b_1b_2).$$
This is well-defined again, as is shown by the following computation
\begin{eqnarray*}
(a_1z_1\otimes b_1)\cdot (a_2\otimes z_2b_2)&=&(-1)^{|b_1||a_2|}(a_1z_1a_2\otimes b_1z_2b_2)\\
&=&(-1)^{|b_1||a_2|+|z_1||a_2|+|b_1||z_2|}(a_1a_2z_1\otimes z_2b_1b_2)\\
&=&(-1)^{|b_1||a_2|+|z_1||a_2|+|b_1||z_2|}(a_1a_2\otimes z_1z_2b_1b_2)\\
&=&(-1)^{|b_1||a_2|+|z_1||a_2|+|b_1||z_2|+|z_1||z_2|}(a_1a_2\otimes z_2z_1b_1b_2)\\
&=&(-1)^{|b_1||a_2|+|z_1||a_2|+|b_1||z_2|+|z_1||z_2|}(a_1a_2z_2\otimes z_1b_1b_2)\\
&=&(-1)^{(|b_1|+|z_1|)(|a_2|+|z_2|)}(a_1a_2z_2\otimes z_1b_1b_2)\\
&=&(a_1\otimes z_1b_1)\cdot (a_2z_2\otimes b_2)
\end{eqnarray*}
for $a_1,a_2\in A$, $b_1,b_2\in B$, $z_1,z_2\in Z$ all homogeneous.
Then, 
$$d_{A\otimes B}=d_A\otimes \textup{id}_B+\id_A\otimes d_B$$
is well-defined again. Indeed, since $d_A(z)=d_B(z)$ for all homogeneous $z\in Z$, we get 
\begin{eqnarray*}
d_{A\otimes B}(az\otimes b)&=&d_A(az)\otimes b+(-1)^{|az|}az\otimes d_B(b)\\
&=&d_A(a)z\otimes b+(-1)^{|a|}ad_A(z)\otimes b+(-1)^{|az|}az\otimes d_B(b)\\
&=&d_A(a)\otimes zb+(-1)^{|a|}a\otimes d_A(z)b+(-1)^{|az|}a\otimes zd_B(b)\\
&=&d_A(a)\otimes zb+(-1)^{|a|}(a\otimes d_B(z)b+(-1)^{|z|}a\otimes zd_B(b))\\
&=&d_A(a)\otimes zb+(-1)^{|a|}(a\otimes d_B(zb))\\
&=&d_{A\otimes B}(a\otimes zb)
\end{eqnarray*}
for all homogeneous $a\in A$ and $b\in B$. This shows that $d_{A\otimes_ZB}$ is well-defined. 
The fact that $d_{A\otimes_ZB}^2=0$ is trivial, and actually follows from the classical case,
such as the fact that it verifies Leibniz' rule. 
\dickebox

\begin{Lemma} \label{ZtensorZisingradedcenter}
Let $(A,d_A)$ and $(B,d_B)$ be two dg-simple dg-algebras with 
$Z_A=Z_{gr}(A,d_A)$ and $Z_B=Z_{gr}(B,d_B)$. Let $Z$ be a common dg-subsalgebra of $Z_A$ and $Z_B$. 
Then $Z_A\otimes_ZZ_B$ is in the graded 
centre of $(A\otimes_ZB,d_{A\otimes_ZB})$. 
\end{Lemma}

Proof.
Indeed, $Z_A\otimes_ZZ_B$ is in the graded center, as 
\begin{eqnarray*}
(z_1\otimes z_2)\cdot (a\otimes b)&=&(-1)^{|a||z_2|}(z_1a\otimes z_2b)\\
&=&(-1)^{|a||z_2|+|z_1||a|+|z_2||b|}(az_1\otimes bz_2)\\
&=&(-1)^{|a||z_2|+|z_1||a|+|z_2||b|+|b||z_1|}(a\otimes b)(z_1\otimes z_2)\\
&=&(-1)^{|z_1\otimes z_2||a\otimes b|}(a\otimes b)(z_1\otimes z_2)
\end{eqnarray*}
This shows the lemma. \dickebox

\section{Remarks on acyclic dg-algebras}

\label{acyclicsection}

\begin{Lemma} \label{acyclicityispreservedunderhomos}
Let $(A,d)$ and $(B,\partial)$ be dg-algebras and suppose that there 
is a dg-algebra homomorphism $\varphi:(B,\partial)\lra (A,d)$. 
If $(B,\partial)$ is acyclic, then $(A,d)$ is acyclic as well. 
\end{Lemma}

Proof. Recall from 
Theorem~\ref{acyclicalgebras} that a dg-algebra $(C,\delta)$ is acyclic 
if and only if $1_C\in\im(\delta)$. Let hence $b\in B$ be homogeneous with $\partial(b)=1_B$. 
Then $$1_A=\varphi(1_B)=\varphi(\partial(b))=d(\varphi(b))\in\im(d)$$
and hence $(A,d)$ is acyclic as well. (We may also argue that $(A,d)$ becomes a 
left $(B,\partial)$-dg-module, and Theorem~\ref{acyclicalgebras} applies.)
\dickebox

\begin{Cor}\label{acyclicgradecentreisacyclic}
Let $(Z,\delta)$ be an acyclic dg-field. Then any dg-algebra $(A,d)$ containing $(Z,\delta)$ 
as a subalgebra is acyclic as well. 
\end{Cor}

Indeed, Lemma~\ref{acyclicityispreservedunderhomos} applies immediately. \dickebox

\begin{Rem} \label{acyclicDelta}
Let $(A,d)$ be a dg-algebra and let $(S_A,\delta_A)$ be a dg-simple dg-module over $(A,d)$.
Then by Proposition~\ref{endosofsimpleisdgfield} the dg-algebra
$$(\End_A^\bullet(S_A,\delta_A),d_\Hom)=:(D_A,\partial_A)$$
is a dg-division algebra, and by Theorem~\ref{dgdivisionalgebraclassification}
we get that either $(D_A,\partial_A)$ is acyclic or $\partial_A=0$. 

If $\partial_A=0$, then every graded endomorphism of $S_A$ automatically 
graded-commutes with the differential of $S_A$. 
\end{Rem}

\begin{Lemma} \label{acyclicsimpleisacyclicendo}
Let $(A,d)$ be a dg-algebra and let $(S_A,\delta_A)$ be a dg-simple dg-module over $(A,d)$. 
Then $(S_A,\delta_A)$ is acyclic if the dg-division algebra
 $(\End_A^\bullet(S_A,\delta_A),d_\Hom)$
is acyclic. 
\end{Lemma}

Proof. Put $$(\End_A^\bullet(S_A,\delta_A),d_\Hom)=:(D_A,\partial_A)$$
If  $(D_A,\partial_A)$ is acyclic, we get that $1_{D_A}\in\im(\partial_A)$.
Hence the identity endomorphism of $S_A$ is a cycle. But this is equivalent with the fact 
that the identity is homotopic to $0$ (cf e.g. \cite[Definition 3.4.6]{Yekutielibook}). 
This then shows that actually $(S_A,\delta_A)$ is acyclic as well. 
\dickebox

\begin{Example} \label{twobytwomatrixalgebras}
The converse of Lemma~\ref{acyclicsimpleisacyclicendo} is false. 
Our favorite example 
$$A=\left(\begin{array}{cc}K&K\\K&K\end{array}\right)$$
with differential 
$$d(\left(\begin{array}{cc}x&u\\z&y\end{array}\right))=\left(\begin{array}{cc}z&y-x\\0&z\end{array}\right)$$
is a counterexample. The differential is non zero, obviously. The graded centre is
just the scalar multiples of the identity.

Also, there is a unique dg-simple dg-left module $(S,\delta)$ over $(A,d)$, namely the right matrix column. 
The differential on this dg-simple is non zero. Actually, $(S,\delta)$ is acyclic. 
The endomorphism ring of this dg-simple $(S,\delta)$ is just $K$, the graded centre.
However, the differential on the graded centre is $0$. 
\end{Example}

\begin{Rem} There is an alternative proof of Lemma~\ref{acyclicsimpleisacyclicendo}. 
Suppose that  $(D_A,\partial_A)$ is acyclic. Since 
$(S_A,\delta_A)$ is a dg-right module over $(D_A,\partial_A)$ 
(cf \cite[Lemma 22.13.3]{stacksproject} or \cite[Proposition 3.3.17]{Yekutielibook}), and since by 
Theorem~\ref{acyclicalgebras} we get that all dg-modules over acyclic dg-algebras are acyclic, also
$(S_A,\delta_A)$ is acyclic. 
\end{Rem}

\begin{Rem}
If we had that any dg-module over a dg-division ring admits a basis in the cycles, then we also had that
any dg-module over a dg-division ring is either acyclic or has differential $0$. Counterexamples
are obiquitous.  

Our favorite example Example~\ref{twobytwomatrixalgebras}
is a counterexample. The differential is non zero, obviously. The graded centre is
just the scalar multiples of the identity. If the statement was correct, then 
the algebra would be isomorphic to the matrix ring with zero differential, which is obviously false. 
\end{Rem}

\begin{Lemma}\label{acyclictensoracyclic}
Let $(A,d_A)$ and $(B,d_B)$ be  dg-algebras with $Z_{gr}(A,d_A)\supseteq Z$ and $Z_{gr}(B,d_B)\supseteq Z$
for some graded commutative dg-algebras $(Z,\partial)$. 
Suppose that $(A,d_A)$ or $(B,d_B)$ is acyclic.
Then $(A\otimes_ZB,d_{A\otimes_ZB})$ is acyclic as well.  
\end{Lemma}

Proof. Indeed, by Lemma~\ref{characterizeacylic} we need to show that 
$1_{A\otimes_ZB}\in\im(d_{A\otimes_ZB})$. 
Since $(A,d_A)$ or $(B,d_B)$ is acyclic, there is a homogeneous 
element $z_A\in A$ or $z_B\in B$ with 
$d_A(z_A)=1_A$ or $d_B(z_B)=1_B$.   
However, $d_{A\otimes_ZB}=d_A\otimes_Z \textup{id}_B+\textup{id}_A\otimes_Zd_B$. 
Suppose that $(A,d_A)$ is acyclic. 
But then 
$$d_{A\otimes_ZB}(z_A\otimes 1_B)=d_A(z_A)\otimes 1_B=1_A\otimes_Z1_B$$
This shows the statement. \dickebox

\section{A differential graded Jacobson-Chevalley density theorem}

\label{dgJacobsonsect}

A very basic result in ring theory is the Chevalley-Jacobson density theorem. 
A graded version of the Chevalley-Jacobson density theorem 
can be found in Chen et al. \cite{Chen}. The result actually follows from earlier work
of Liu, Beattie and Fang~\cite{Beattie}.

\begin{Theorem} \cite{Chen, Beattie} \label{gradeddensitytheorem}
Let $A$ be a group graded algebra. 
Let $M$ be a gr-simple graded $A$-module and let $D=\End_{A-graded}(M)$. 
Let $x_1,\dots,x_k$ be homogeneous $D$-independent elements of $M$, and let $y_1,\dots,y_k$
be any elements of $M$. Then there is $a\in A$ with $ax_i=y_i$ for any $i\in\{1,\dots,k\}$.  
\end{Theorem}

We shall need to show that if $(A,d)$ is a dg-algebra 
and $(M,\delta)$ is a dg-simple dg-module over $(A,d)$, then
$\ker(\delta)$ is a $\Z$-gr-simple $\ker(d)$-module. 
In the special case of an acyclic $(A,d)$-algebra, this follows from 
Theorem~\ref{acyclicalgebras}, due to Aldrich and Garcia-Rozas~\cite{Tempest-Garcia-Rochas}.


As a consequence, if $(A,d)$ is acyclic, then a dg-module $(M,\delta)$ over $(A,d)$ 
is dg-simple if and only if $\ker(\delta)$ is a graded-simple $\ker(d)$-module. 
We hence may apply Theorem~\ref{gradeddensitytheorem} to dg-simple dg-modules
in this case.

\begin{Rem}
\begin{itemize}
\item
An unpublished result \cite{Bergman2} due to G.M.Bergman shows that for $\Z$-graded rings, the Jacobson radical 
(ungraded version) is homogeneous. This means that any simple module is automatically graded. 
\item
Recall that Bahturin, Zaicev and Sehgal classified in \cite{BahturinZaicevSehgal}
finite-dimensional simple $G$-graded $K$-algebras $A$ for a group $G$ and an algebraically closed field $K$, 
subject to some hypotheses with respect to  $G$ and to
the base field. In particular, if either $K$ is of characteristic $0$ or 
the order of any finite subgroup of $G$ is coprime to the characteristic of $K$, then 
$A$ is a matrix algebra over a graded skew-field $K^\alpha H$ for some finite 
subgroup $H$ of $G$ and $\alpha$ a $2$-cocycle with values in $K^\times$. 

If $G=\Z$, then there is no non trivial finite subgroup. Hence, considering finite dimensional graded 
simple algebras we are left with gradings on full matrix algebras over $F$.  
\item
Note that $\ker(d_{A\otimes_ZB})$ is in general strictly bigger that 
$\ker(d_A)\otimes_ZB+A\otimes_Z\ker(d_B)$. Hence, the equivalence of categories in 
Theorem~\ref{acyclicalgebras} does not behave well with respect to tensor products. 
\end{itemize}
\end{Rem}

\begin{Theorem}\label{dgdensitytheorem}
Let $(A,d)$ be an acyclic dg-algebra, i.e. a dg-algebra with $H(A,d)=0$.
Let $(M,\delta)$ be a dg-simple dg-module over $(A,d)$ and let 
$$(D,\partial):=\End_A^\bullet((M,\delta),d_\Hom).$$
Then $(D,\partial)$ is a dg-division algebra, and also a $\Z$-gr-division algebra.
Moreover, for each family $x_1,\dots,x_k$ 
of $D$-independent elements of $\ker(\delta)$  and each family 
$y_1,\dots,y_k$ of elements of $\ker(\delta)$, there is an element $a\in\ker(d)$
with $ax_i=y_i$ for all $i\in\{1,\dots,k\}$.
\end{Theorem}

Proof. Let $N:=\ker(\delta)$ and $B:=\ker(d)$ to shorten the notation. 
As $(M,\delta)$ is dg-simple over $(A,d)$, by the equivalence of categories in 
Theorem~\ref{acyclicalgebras}.(\ref{acyclicalgebras5})
we get that $N$ is $\Z$-graded simple as graded $B$-module.
Further, $\End_{B-graded}(N)\simeq D=\End_A^\bullet((M,\delta),d_\Hom)$ again 
by the equivalence of categories in
Theorem~\ref{acyclicalgebras}.(\ref{acyclicalgebras5}).
Hence, $D$ is also a $\Z$-gr-division algebra. 
The statement now follows directly from Theorem~\ref{gradeddensitytheorem}. \dickebox


\begin{thebibliography}{88}


\bibitem{Tempest-Garcia-Rochas}
S. Tempest Aldrich and J. R. Garcia Rozas, {\em Exact and Semisimple Differential Graded Algebras},
Communications in Algebra {\bf 30} (no 3) (2002) 1053-1075.

\bibitem{BahturinZaicevSehgal}
Yu. A. Bahturin, M. V. Zaicev, Sudarshan K. Sehgal, {\em Finite-dimensional simple graded algebras},
(English version) Sbornik Mathematics {\bf 199:7} 965-983.



\bibitem{Bergman2}
George M. Bergman, {\em On Jacobson radicals of graded rings}, preprint (1975) \\
{\tt https://math.berkeley.edu/~gbergman/papers/unpub/J\textunderscore G.pdf}


\bibitem{Cartandg}
Henri Cartan, {\em DGA-alg\`ebres et DGA-modules, } Séminaire Henri Cartan, tome 7, no 1 (1954-1955),
exp. no 2, p. 1-9.

\bibitem{Chen}
Tung-Shyan Chen, Chin-Fang Huang, and Jing-Whei Liang, {\em Extended Jacobson Density Theorem for Graded Rings with Derivations and Automorphisms}, Taiwanese Journal of Mathematics {\bf 14} no 5 (2010) 1993-2014.

\bibitem{CR1}
Charles W. Curtis and Irving Reiner, {\sc Methods of Representation Theory}, Vol 1,
John Wiley Interscience 1981.


\bibitem{Goodbody}
Isambard Goodbody, {\em Reflecting perfection for finite dimensional differential graded algebras}. preprint october 5, 2023; arxiv:2310.02833.


\bibitem{GoodearlStafford}
Ken Goodearl and Tobi Stafford, {\em The graded version of Goldie's theorem}, Contemporary Math. {\bf 259}, (2000) 237-240.




\bibitem{Beattie}
S. X. Liu, Margaret Beattie, and H. J. Fang, {\em Graded division rings and the Jacobson density theorem}, Journal of the Beijing Normal University (Natural Science) {\bf 27} (2) (1991) 129-134.


\bibitem{gradedrings}
Constantin Nastasescu and Fred van Oystaen, {\sc Graded Ring Theory}. North Holland 1982; Amsterdam



\bibitem{Orlov1}
Dimitri Orlov, {\it Finite dimensional differential graded algebras and their geometric realisations}, Advances in Mathematics {\bf 366} (2020) 107096, 33 pp.

\bibitem{Orlov2}
Dimitri Orlov, {\em Smooth DG algebras and twisted tensor product .} arxiv 2305.19799






\bibitem{stacksproject}
The stacks project. {\tt https://stacks.math.columbia.edu/browse}

\bibitem{Yekutielibook}
Amnon Yekutieli, {\sc Derived categories}, Cambridge studies in advanced 
mathematics {\bf 183}, Cambridge university press, Cambridge 2020.


\bibitem{dgorders}
Alexander Zimmermann, {\em Differential graded orders, their class groups and id\`eles}, preprint December 30, 2022; 32 pages

\bibitem{dgBrauer}
Alexander Zimmermann, {\em Differential graded Brauer groups}, preprint March 31, 2023; final version 11 pages; to appear in Revista de la Union Matematica Argentina.

\bibitem{dgGoldie}
Alexander Zimmermann, {\em Ore Localisation for differential graded rings; Towards
Goldie's theorem for differential graded algebras}, Journal of Algebra {\bf 663} (2025) 48-80. 

\end{thebibliography}
\end{document}